\documentclass[12pt,a4paper]{amsart} 
\usepackage{amssymb,latexsym}

\def\N{\mathbb N} 
 
\def\R{\mathbb R} 
\def\D{\mathbb D} 
 
\def\T{\mathbb T} 
 
\def\H2{{\mathcal H}^2(\D)} 
\newtheorem{lem}{LEMMA} 
\newtheorem{theo}[lem]{THEOREM} 
 
\newtheorem{coro}[lem]{COROLLARY}

\begin{document} 
\title[ Multilinear Hankel operators]{Truncation of Multilinear Hankel operators} 
\author{Sandrine Grellier \& Mohammad Kacim} 
\address{MAPMO\\Universit\'e d'Orl\'eans\\ Facult\'e des Sciences\\ D\'epartement de Math\'ematiques\\BP 6759\\ F 45067 ORLEANS C\'edex 2 
\\FRANCE} 
\email{ grellier@labomath.univ-orleans.fr\\ kacim@labomath.univ-orleans.fr } 
\thanks{Authors partially supported by the 2002-2006 IHP Network, Contract Number: HPRN-CT-2002-00273 - HARP\\ 
The authors would like to thank Joaquim Bruna who suggested this 
problem.} 
\begin{abstract} We extend to multilinear Hankel operators the fact that truncation of  bounded Hankel operators is bounded. 
We prove and use a continuity property of a kind of bilinear Hilbert 
transforms on product of Lipschitz spaces and Hardy spaces. 
\end{abstract} 

\maketitle 
\section{Statement of the results} 
 
In this note, we  prove that truncations of bounded multilinear 
Hankel operators are bounded. This extends the same property for 
linear Hankel operators, a result obtained by \cite{BB}, which we 
first recall. A matrix $B=(b_{mn})_{m,n\in\N}$ is called of {\sl 
Hankel type} if $b_{mn}=b_{m+n}$ for some sequence $b\in l^2(\N)$. 
We can identify $B$ with an operator acting on $l^2(\N)$. 
Moreover, if we identify $l^2(\N)$ with the complex Hardy space 
$\H2$ of the unit disc, then $B$ can be realized as the integral 
operator, called {\sl Hankel operator} and denoted by $H_b$, which 
acts on $f\in \H2$ by 
$$H_bf(z)= 
\frac{1}{2\pi}\int_{\T}\frac{b(\zeta)f(\overline{\zeta})}{1-\overline{\zeta}z}d\sigma(\zeta).$$ 
In other words, $H_bf=\mathcal C(b\check f)$ where $\mathcal C$ 
denotes the Cauchy integral, $\check 
f(\zeta):=f(\overline{\zeta})$, $\zeta\in\T$. The {\sl symbol} $b$ 
  of the Hankel operator is given by $b(\zeta):=\sum_{k=0}^\infty 
b_k\zeta^k$. If $f(z)=\sum_{n\in\N}a_nz^n$, one has 
$$H_bf(z)=\sum_{m\in\N}(\sum_{n\in\N} a_nb_{m+n})z^m.$$ Now, we 
consider truncations of matrices defined  as follows. For 
$\beta,\gamma\in\R$, the truncated matrix $\Pi_{\beta,\gamma}(B)$ 
is the matrix whose $(m,n)$ entry  is $b_{mn}$ or zero, depending 
on the fact that  $m\ge \beta n+\gamma$ or not. It is proved in 
\cite{BB} that such truncations, for $\beta\neq -1$, preserve the 
boundedness for Hankel operators. The proof consists in  showing 
that truncations are closely related to  bilinear periodic Hilbert 
transforms. One then uses the theorem of Lacey-Thiele (see 
\cite{LT1}, \cite{LT2}, \cite{LT3}) in the periodic setting. We 
are interested in the same problem  for multilinear Hankel 
operators. For $n\in\N$, we define the multilinear Hankel operator 
$H^{(n)}_b$ as follows. Let $f_1,\dots,f_n \in\H2$, 
\begin{eqnarray*}H_b^{(n)}(f_1,\dots,f_n)(z)&=&\frac{1}{2\pi}
\int_{\T}\frac{b(\zeta)f_1(\overline{\zeta})\dots f_n(\overline{\zeta})} 
{1-\overline{\zeta}z}d\sigma(\zeta)\cr&=&H_b(f_1\times\dots\times 
f_n)(z). 
\end{eqnarray*} 
When equipped with the canonical basis of $\H2$, this operator has 
a matrix $B$ with entries in $\N^{n+1}$, which we denote by 
$B=(b_{i_0,\dots,i_{n}})_{i_0,\dots,i_{n}\in\N}$. We speak of 
$(n+1)$-dimensional infinite matrices (so that a usual matrix is a 
$2$-dimensional matrix in our terminology). Its action on $n$ 
vectors $a^1\dots,a^n$ gives the vector whose 
 $m$-th coordinate is 
$$\sum_{i_1,\dots,i_{n}}b_{m,i_1,\dots,i_n}a^1_{i_1}\dots a^n_{i_{n}}.$$ 
In the case of the operator $H_b^{(n)}$, the matrix $B$ is a 
$(n+1)$-dimensional matrix with entries which are constant on the 
hyperplanes $i_0+\dots+i_{n}=c$. Such a matrix is called a 
$(n+1)$-dimensional Hankel matrix. 
 
We consider truncations of $(n+1)$- dimensional matrices as 
follows. For ${\bf\beta}\in \R^n$ and $\gamma\in\R$, 
${\bf\beta}=(\beta_1,\dots,\beta_n)$, $\Pi_{{\bf 
\beta},\gamma}(H_b^{(n)})$ denotes the $(n+1)$-dimensional matrix 
whose $(i_0,\dots,i_{n})$ entry is $b_{i_0+\dots+i_{n}}$  if 
$\beta_1 i_1+\dots+\beta_n i_n+\gamma\le i_{0}$ and zero 
otherwise. In this note, we consider the simplest case where 
$\beta=(1,\dots,1)$ and $\gamma=0$ which is denoted by $\Pi_{1,0}$. We will study the general case in a
foregoing paper. Our main result is the following. 
 
\begin{theo}\label{main} 
If $H_b^{(n)}$ is a bounded multilinear Hankel operator from 
$\left(\H2\right)^n$ into $\H2$ then so is its truncated operators 
$\Pi_{1,0}(H_b^{(n)})$. 
\end{theo} 
 
Theorem \ref{main} is  a corollary of an estimate on a kind of bilinear 
Hilbert transform in the periodic setting which is of independent 
interest. Let us first give some notations. The usual Lipschitz 
spaces of order $\alpha$ are denoted by $\Lambda_\alpha(\T)$, 
while  $\mathcal H^p(\T)$ denotes the real Hardy space, $p>0$. 
 
Let us finally recall that for $f$ and $b$ trigonometric 
polynomials on the torus, the periodic bilinear Hilbert transform 
of $f$ and $b$ is given by 
$${\mathcal H} (b,f)(x)= p.v.\int_{\T} b(x+t)f(2t) {dt \over \tan{x-t \over 2}}.$$ 
Lacey-Thiele's Theorem, once transferred to the periodic setting, 
is  the following. 
 
\begin{theo}{\bf \cite{BB}}\label{BB} 
Let $1<p,q\le \infty$ with $\displaystyle\frac 1r=\frac 1p+\frac 1q<\frac 32$. Then, there exists a constant $C>0$ so that, 
for any trigonometric polynomials $f$ and $b$, 
$$\Vert {\mathcal H} (b,f)\Vert _r\le C\Vert f\Vert _p\Vert b\Vert_q.$$ 
\end{theo} 
 
We adapt the definition to our setting, and define 
\begin{equation}\label{bilinear2} 
  \widetilde{\mathcal  H}(b,f)(x)= \int_{\T} \left (b(x+t)-b(2x)\right)f(2t) \frac{dt}{ \tan{\frac{x-t} 2}}. 
\end{equation}

We prove the following. 
\begin{theo}\label{bilinear-p} 
Let $1<p<\infty$, $0<q<p$ and  $\alpha=\frac 1q-\frac 1p$. There exists a constant $C>0$ so that, for any sufficiently smooth functions $b\in\Lambda_\alpha(\T)$ and 
$f\in\mathcal H^p(\T)$ 
\begin{equation}\label{bilin} 
\Vert \widetilde{\mathcal  H} (b,f)\Vert _{\mathcal H^p}\le 
C\Vert f\Vert _{\mathcal H^q}\Vert b\Vert_{\Lambda_\alpha}. 
\end{equation} 
\end{theo} 
 
We remark that the limiting case $b\in L^\infty (\T)$ is given by 
the Lacey-Thiele theorem, Theorem \ref{BB}. 
 
\medskip 
 
Let us come back to holomorphic functions and to truncations. 
Denote by $\Lambda_\alpha(\D)$, $\alpha>0$, the space of functions 
which are holomorphic in $\D$ and whose boundary values are in 
$\Lambda_\alpha(\T)$. Denote also by $\mathcal H^p(\D)$ the 
complex Hardy space on the unit disc, $p>0$. Recall that the dual 
of $\mathcal H^p(\D)$ is $\Lambda_\alpha(\D)$, for $\displaystyle 
p=\frac{1}{\alpha+1}$ (\cite{D}). As an easy consequence of 
duality and factorization, one obtains that the Hankel operator 
$H_b$ is bounded from $\mathcal H^q(\D)$ into $\mathcal H^p(\D)$, 
with $q<p$ and $p>1$, if and only if the symbol $b$ is in 
$\Lambda_\alpha(\D)$ with $\alpha=\frac 1q-\frac 1p$. More precisely, there exists a constant $C$ such that, 
for all holomorphic polynomials $f$, 
\begin{equation}\label{hankel} 
\|H_b(f)\|_p\leq C\|b\|_{\Lambda_\alpha(\D)}\times \|f\|_{\mathcal 
H^{q}(\D)}. 
\end{equation}

Theorem \ref{bilinear-p} has the following corollary, which gives 
the link with truncations.

\begin{coro} 
Let $1<p<\infty$ and $0<q<p$. Let  $b\in\Lambda_\alpha(\D)$. Then the operator 
 $\Pi_{1,0}(H_b)$ is  bounded  from $\mathcal H^q(\D)$ into $\mathcal H^p(\D)$. 
\end{coro} 
 
So, if $H_b$ is a bounded operator, its truncate 
$\Pi_{1,0}(H_b)$ is also bounded. 
 
Let us deduce Theorem \ref{main} from the corollary. It is clear 
that $$H_b^{(n)}(f_1,\dots,f_n)(z)=H_b(f_1\times\dots\times 
f_n)(z).$$ Using the factorisation of functions in Hardy classes, 
we know that $H_b^{(n)}$ is bounded as an operator from 
$\left(\H2\right)^n$ into $\H2$ if and only if the Hankel operator 
$H_b$ is bounded  from $\mathcal H^{2/n}(\D)$ into $\H2$, that is, 
if and only if $b$ is in $\Lambda_\alpha(\D)$ for $\alpha=\frac 
{n-1}2$ . To conclude, we use the fact that the truncation 
$\Pi_{1,0}$ of $H_b^{(n)}$ corresponds to the truncation 
$\Pi_{1,0}$ of $H_b$, as it can be easily verified. 
 
The remainder of the paper is organized as follows. In the next 
section, we deduce the corollary from Theorem \ref{bilinear-p}. In 
the last one, we prove Theorem\ref{bilinear-p}. Let us emphasize the 
fact that this last proof does not use Lacey-Thiele Theorem, and 
is elementary compared to it.

\section {The link between truncations and bilinear Hilbert 
transforms} 
 
We prove the corollary. It is sufficient to prove that,  
for $b$ and $f$ trigonometric polynomials, 
\begin{equation}\label{ineq} 
\|\Pi_{1,0}H_b(f)\|_p\leq 
C\|b\|_{\Lambda_\alpha(\D)}\times \|f\|_{\mathcal H^{q}(\D)} 
\end{equation} 
for some constant $C$ which is independent of $b$ and $f$.  
 
Let $b$ and $f$ be two triginometric polynomials. 
Assume that $f(z)=\sum_{n\in\N} a_n z^n$ and denote by $F$ the 
function defined on the torus by $F(x)=f(e^{-ix})$. It is elementary to see that $F$ and $f$ have the 
same norm in $\mathcal H^{q}(\T)$. Moreover, an elementary 
computation (which is already in \cite{BB}) shows that the analytic 
part of $\mathcal{H}(b,F)(x)$ is equal to 
$$\sum_{n\in\N}\sum_{m\in\N} a_nb_{m+n} sign((m-n)) 
e^{i2mx}.$$ So it is sufficient to prove that 
$\mathcal{C}(\mathcal{H}(b,F))$ satisfies the 
desired estimate, 
\begin{equation}\label{gamma} 
\|\mathcal{C}(\mathcal{H}(b,F))\| _p\leq 
C\|b\|_{\Lambda_\alpha(\D)}\times \|f\|_{\mathcal H^{q}(\D)}. 
\end{equation} 
 
We want to replace $\mathcal{H}(b, F)$ by $\widetilde{\mathcal{H} }(b,F)$, for which we have such 
an estimate given in Theorem \ref{bilinear-p}. Let us look at the 
difference, which is given, up to a constant, by the Cauchy 
projection of 
$$x\mapsto b(2x)\int_{\T} F(2t) \frac{dt}{ \tan{\frac{x-t} 2}}=b(2x)F(2x)$$ 
since $f$ has only non zero coefficients for positive frequencies. We 
recognize $H_b(f)(z^2)$, whose norm in $\mathcal H^p(\D)$ 
coincides with the one of $H_b(f)$. To conclude for (\ref{gamma}), 
we use (\ref {bilin}) and (\ref {hankel}). 
 

%
%
\section{Proof of the Theorem 3.} 
  
When $q>1$, then the kernel of $\widetilde{\mathcal H}$ is bounded , up to a constant $c\Vert b\Vert_\alpha$, by the Riesz 
potential $|x-y|^{-1+\alpha}$, and the estimate follows directly. 
Let us now concentrate on $q\leq 1$, for which we can use the 
atomic decomposition. By the atomic decomposition Theorem of 
$H^q(\T)$, it suffices to consider the action of $\widetilde{\mathcal H}(b,.)$ on $H^q(\T)$-atoms. 
Let $a$ be a $H^q(\T)$-atom. 
 
If $a$ is the constant atom, or if $a$ is a non constant atom 
which is supported in some interval $I$ of the torus of length 
bigger than $\frac \pi4$, then, for all $r>1$, its $L^r$ norm is 
uniformly bounded. It follows at once that the $L^p$ norm of 
$\widetilde{\mathcal H}(b,a)$ is also uniformly bounded. 
 
We assume now that $a:t\mapsto \tilde a(2t)$ is an atom supported in some interval $I$ 
on the torus, centered at $x_I$ and of radius $r<\frac {\pi}4$. 
Denote by $\tilde I$ the interval centered at $x_I$ and of radius 
$2r$, and by $2\tilde I$ the ball centered at $2x_I$ and of radius 
$4r$. We first consider the case when $0<\alpha<1$. We write
\begin{eqnarray*} 
\widetilde{\mathcal H}(b,\tilde a)(x) &=&  \underbrace{ p.v.\int_{\T} (b(x+t)-b(2x))a(t) 
{dt \over \tan{x-t \over 2}}1\!\text{l}_{x\in \tilde I}}_{A_1(x)}\\ 
&+&   \underbrace{p.v.\int_{\T} (b(x+t)-b(2x))a(t) {dt \over \tan{x-t \over 2}}
1\!\text{l}_{x\in {\tilde I}^c}}_{A_2(x)} 
\end{eqnarray*} 
 
We prove that both $A_1$ et $A_2$ are $L^p$-functions. To 
prove that  $A_1 \in L^p$, we write 
$\vert A_1(x)\vert\le \Vert b\Vert_\alpha \mathcal I_\alpha(|a|)(x)1\!
\text{l}_{x\in \tilde I}$ where $\mathcal I_\alpha$ denotes the 
fractional integral
related to the Riesz potential $|x-y|^{-1+\alpha}$. So, by Minkowski inequality,
$$\Vert A_1\Vert_{L^p}\le c \Vert b\Vert_\alpha\Vert \mathcal
I_\alpha(a)||_{L^s}\times |I|^{1/p-1/s}$$ for any $s>p$. 
We choose $s$ large enough to have $\alpha+1/s<1$ and $r>1$ so that $1/r=\alpha+1/s$. For these choices, we get
$$\Vert A_1\Vert_{L^p}\le c \Vert b\Vert_\alpha\Vert a\Vert_{L^r}\times |I|^{1/p-1/s}\le c \Vert b\Vert_\alpha.$$

It remains to consider the term denoted by $A_2$. For this term, we write 
$$b(x+t)-b(2x)=[b(x+t)-b(x+x_I)]+[b(x+x_I)-b(2x)].$$ 
The corresponding terms are denoted by $A_2^{(1)}$ and $A_2^{(2)}$ respectively. 
For the first term $A_2^{(1)}$, we use that 
$$|b(x+t)-b(x+x_I)|\le c\Vert b\Vert_\alpha |I|^\alpha$$ for any $t\in I$ and that 
$\vert \tan{x-t \over 2}\vert \ge  C|x-x_I|$ when $t\in I$ and $x\in\tilde I^c$. So, 
\begin{eqnarray*} 
\vert A_2^{(1)}(x)\vert&\le&c  \Vert b\Vert_\alpha |I|^\alpha
\times \frac{1}{|x-x_I|}\left(\int_{I} |a(t)|dt\right) 1\!\text{l}_{\tilde
I^c}(x)\\ 
&\le & c\Vert b\Vert_\alpha |I|^{\alpha+1-1/q}\frac{1}{|x-x_I|}1\!
\text{l}_{\tilde I^c}(x). 
\end{eqnarray*} 
Taking the $L^p$-norm, it gives, as $p>1$, 
$$\Vert A_2^{(1)}\Vert_{L^p}\le c\Vert b\Vert_\alpha |I|^{\alpha+1-1/q} 
\left(\int_{|x-x_I|\ge 2r}\frac{1}{|x-x_I|^p}dx\right)^{1/p} 
\le c \Vert b\Vert_\alpha.$$ 
For the second part, we use the fact that $a$ has vanishing moment of order $m:=\left[\frac1q\right]-1$ so that one can substract to 
$t\to {1\over \tan{x-t \over 2}}$ its Taylor expansion of order $m$ at point $x_I$ without changing the value of 
$A_2^{(2)}$. As the corresponding difference is bounded by $\displaystyle \frac{|I|^{m+1}}{|x-x_I|^{m+2}}$ for $t\in I$ and 
$x\in\tilde I^c$, it allows to obtain 
\begin{eqnarray*} 
\vert A_2^{(2)}(x)\vert &\le & c\Vert b\Vert_\alpha  |x-x_I|^\alpha\times 
\frac{|I|^{m+1}}{ |x-x_I|^{m+2}} \times 
\left(\int_{\tilde I}|a|\right)1\!\text{l}_{\tilde I^c}(x)\\ 
&\le &c\Vert b\Vert_\alpha  |x-x_I|^{\alpha-m-2}\times |I|^{m+2-1/q}1\!\text{l}_{\tilde I^c}(x). 
\end{eqnarray*} 
Eventually, it gives 
$$\Vert A_2^{(2)}\Vert_p \le c\Vert b\Vert_\alpha |I|^{m+2-1/q}\times \left(\int_{|x-x_I|\ge 2r}|x-x_I|^{(\alpha-m-2)p}dx\right)^{1/p}.$$ 
This last integral is convergent since, as $m=\left[\frac1q\right]-1$ and $\alpha=\frac 1q-\frac 1p$, $(\alpha-m-2)p=-1+(\frac 1q-m-2)p<-1$. 
So, we obtain
$$\Vert A_2^{(2)}\Vert_p \le c\Vert b\Vert_\alpha.$$ 
So, we have proved that $\Vert \widetilde{\mathcal H}(b,a)\Vert_{L^p}\le  c\Vert b\Vert_\alpha$ for any $H^q$-atom $a$. 
It proves that $\widetilde{\mathcal H}(b,\cdot)$ maps $H^q(\T)$ 
into $L^p(\T)$ boundedly. It ends the proof of the theorem in the case $0<\alpha<1$. 
 
Now, we illustrate the method for larger values of $\alpha$ by considering the case $1\le \alpha<2$. We write 
\begin{eqnarray*} 
\tilde{\mathcal H}(b,a)(x)&=&\int_\T (b(x+t)-b(2x)-(x-t)b'(2x))a(t) \frac{dt}{\tan\frac{x-t}{2}}\cr
&+&b'(2x)K* a(x)\cr
&=& H_1(x)+H_2(x). 
\end{eqnarray*} 
Here $K$ is the $\mathcal C^\infty$-kernel defined by $K(x):=\frac{x}{\tan\frac{x}{2}}$, $x\in\T$. 
The corresponding term $H_2$ is hence in $H^p(\T)$ since, as $K$ is a $\mathcal 
C^\infty$-kernel, $K*a$ is a smooth function (even if $a$ is only a distribution in $H^q(\T)$). 
In particular it belongs to $L^p(\T)$ and so is for $b'(2.)K*a$ with 
$$\Vert b'(2.)K* a\Vert_{H^p}\le C\Vert b'\Vert_\infty\le C\Vert b\Vert_\alpha.$$ 
So, the problem reduces to show that $H_1$ belongs to $L^p(\T)$. 
We write as before $H_1(x)= A_1(x)+A_2(x)$ where $A_1(x)=H_1(x)1\!\text{l}_{\tilde I}(x)$. 
 
To prove that  $A_1 \in L^p(\T)$, 
we write, for $x\in\tilde I$
\begin{eqnarray*} 
(*)&:=&b(x+t)-b(2x)-(t-x)b'(2x)\cr
&=&[b(x+t)-b(2x_I)-(x+t-2x_I)b'(2x_I)]\cr 
&-&[b(2x)-b(2x_I)-(2x-2x_I)b'(2x_I)]\cr 
&+&(t-x)[b'(2x_I)-b'(2x)]\cr 
&:=& \tilde b(x+t)-\tilde b(2x)+(t-x)[b'(2x_I)-b'(2x)] 
\end{eqnarray*} 
where $\tilde b(s):=[b(s)-b(2x_I)-(s-2x_I)b'(2x_I)]\Psi(s)$ with $\Psi$ a smooth function supported in twice of $2\tilde I$ identically $1$
in $2\tilde I$.
We write $A_1=A_1^{(1)}+A_1^{(2)}$. 
For the first term, we remark that $\tilde b$ belongs to $\Lambda_\beta(\T)$ for any $0<\beta<1$ with 
$$\Vert \tilde b\Vert_{\Lambda_\beta}\le c \Vert b\Vert_\alpha|I|^{\alpha-\beta}.$$ This follows from the fact that $\Vert \tilde b'\Vert_\infty\le
c\Vert b\Vert_\alpha|I|^{\alpha-1}$ (since, by the choice of $\Psi$, $\Vert \Psi'\Vert_\infty \le  c|I|^{-1}$).
From the first part of the proof, we get that, 
for $\displaystyle \beta=1-\frac 1p$, the corresponding operator maps $H^1(\T)$ into 
$H^p(\T)$ with
$$\Vert A_1^{(1)}\Vert_{L^p}\le c\Vert \tilde b\Vert_{\Lambda_\beta} \Vert a\Vert_{L^1}\le c|I|^{\alpha-\beta+1-1/q}\le c.$$

For $A_1^{(2)}$, we have $A_1^{(2)}(x)=[b'(2x_I)-b'(2x)]K*a(x)1\!\text{l}_{\tilde I}(x)$. 
So, it gives 
$$\Vert A_1^{(2)}\Vert_{L^p}\le c\Vert 
b\Vert_\alpha |I|^{\alpha-1+1/p}\times\int_\T|a|\le c\Vert b\Vert_\alpha.$$ 
 
To deal with $A_2$, we write
\begin{eqnarray*} 
(*)&=&b(x+t)-b(2x)-(t-x)b'(2x)\cr
&=&[b(x+t)-b(x+x_I)-(t-x_I)b'(x+x_I)]\cr 
&+&[b(x+x_I)-b(2x)-(x_I-x)b'(x+x_I)]\cr 
&+&(t-x)[b'(x+x_I)-b'(2x)]. 
\end{eqnarray*} 
The corresponding terms are denoted by $A_2^{(1)}$, $A_2^{(2)}$ and $A_2^{(3)}$ respectively. 
For the first term $A_2^{(1)}$, we use that $|b(x+t)-b(x+x_I)-(t-x_I)b'(x+x_I)|\le c\Vert b\Vert_\alpha|I|^\alpha$ for any 
$t\in I$ and that $\vert \tan{x-t \over 2}\vert \ge  C|x-x_I|$ when $t\in I$ and $x\in\tilde I^c$. So, the estimate of this term is as 
before. 
For the second part, we use the fact that $a$ has vanishing moment of order $m:=\left[\frac1q\right]-1$ so that one can substract to 
$t\to {1\over \tan{x-t \over 2}}$ its Taylor expansion of order $m$ at point $x_I$ without changing the value of 
$A_2^{(2)}$. As the corresponding difference is bounded by $\displaystyle \frac{|I|^{m+1}}{|x-x_I|^{m+2}}$ for $t\in I$ and 
$x\in \tilde I^c$, it allows to obtain the same estimate as before. 
 
We just have to consider the third term $A_2^{(3)}$. Here, we write that $t-x=[t-x_I]+[x_I-x]$ so that it gives two different 
terms to estimate. In the first, we use again that $a$ has vanishing moment of order less than 
$m$ so that one can substract to $t\to {1\over \tan{x-t \over 2}}$ its Taylor expansion of order $m-1$ at point $x_I$ without 
changing the value of the integral $\displaystyle \int_\T [t-x_I]a(t) {dt \over \tan{x-t \over 2}}.$ So, it gives that the corresponding 
term is bounded by $\Vert b\Vert _\alpha |x-x_I|^{\alpha-1-m-1} |I|^{2-1/q+m}$ and its corresponding $L^p$-norm is bounded by 
$C\Vert b\Vert _\alpha.$ 
 For the very last term, one can substract to $t\to {1\over \tan{x-t \over 2}}$ its
 Taylor expansion of order $m$ at point 
$x_I$ without changing the value of the integral $\displaystyle \int_\T a(t) {dt \over \tan{x-t \over 2}}$, 
it gives that the corresponding term has the same bound as the preceding one. It 
finishes the proof.


\begin{thebibliography}{MTW1} 
\bibitem[BB]{BB}  Bonami, Aline; Bruna, Joaquim {\em On truncations of Hankel and Toeplitz operators.} 
Publ. Mat. 43 (1999), no. 1, 235--250 
 
\bibitem[BL]{BL} Bergh, J\"oran; L\"ofstr\"om, J\"orgen {\em Interpolation spaces. An introduction.} 
Grundlehren der Mathematischen Wissenschaften, No. 223. Springer-Verlag, Berlin-New York, 1976. x+207 pp 
 
\bibitem[CM]{CM} Coifman, R.R. and Meyer, Y.  {\em On commutators of singular integrals and bilinear singular 
integrals,} Trans. Am. Math. Soc., 212, 
315--331 (1975). 
 
\bibitem [D]{D} Duren, Peter L. {\em Theory of $H\sp{p}$ spaces.} Pure and Applied Mathematics, 
Vol. 38 Academic Press, New York-London 1970 
xii+258 

\bibitem[GR]{GR}  Garc\'{\i}a-Cuerva, Jos\'e; Rubio de Francia, Jos\'e L. {\em Weighted norm inequalities and related topics.} 
North-Holland Mathematics Studies, 116. Notas de Matemática [Mathematical Notes], 104. North-Holland Publishing Co., Amsterdam, 1985. x+604 pp. 

\bibitem [GN]{GN} Gilbert, J.E. and Nahmod, A.R.  {\em 
$L^p$-boundedness of time-frequency paraproducts, II,} J. of 
Fourier Anal. Appl., 8 (2), 109-172 (2002). 
 
\bibitem[GT]{GT} Grafakos, L. and Torres, R.H.{\em Multilinear Calder\'{o}n-Zygmund theory,} Adv. in 
Math., 165, 124-164. 
 
\bibitem[KS]{KS} Kenig, C.E. and Stein, E.M. (1999). {\em Multilinear estimates and fractional integration,} 
Math. Research Lett., 6, 1--15 (1999). 
 
\bibitem[L]{L} Lacey, M. (2000). {\em The bilinear maximal function maps in $L^p$ for $2/3 p \leq 1$,} Ann. of Math., 
151(2), 35--57 (2000). 
 
\bibitem [LT1]{LT1} Lacey, M. and Thiele, C. (1997). {\em $L^p$ estimates on the bilinear Hilbert transform, $2< p< \infty$,} 
Ann. of Math., 146, 683--724 (1997). 
 
\bibitem [LT2]{LT2} Lacey, Michael T.; Thiele, Christoph M. 
{\em On Calder\'on's conjecture for the bilinear Hilbert transform. }Proc. Natl. Acad. Sci. USA 95 (1998), no. 9, 4828--4830 
 
\bibitem [LT3]{LT3} Lacey, M. and Thiele, C. {\em On Calder\'on's conjecture,} Annals of Math., 149, 475--496 (1999). 
 
 
\bibitem [MTT]{MTT} Muscalu, C., Tao, T., and Thiele, C. {\em Multilinear operators given by singular multipliers,} J. Am. Math. Soc., 15, 469-496 (2002). 
 
\bibitem [S]{S} Stein, E.M. {\em Harmonic Analysis: Real Variable Methods, Orthogonality, and Oscillatory Integrals,} Princeton University Press, 
Princeton, NJ (1993). 
 
\end{thebibliography}
\end{document}